# Development and Modelling of High-Efficiency Computing Structure for Digital Signal Processing

Annapurna Sharma[1], Hakimjon Zaynidinov[2], Hoon Jae Lee [3]

[1] *Graduate School of Design & IT, Dongseo University, Korea*
[2]*Division of Design & IT, Dongseo University, Korea*
[3]*Division of Computer and Information Engineering, Dongseo University*
*Email*: sharmaannapurna@gmail.com[1] tet2001@rambler.ru[2] , hjlee@dongseo.ac.kr[3] ,

*Abstract*— **The paper is devoted to problem of spline approximation. A new method of nodes location for curves and surfaces computer construction by means of B-splines and results of simulink-modeling is presented. The advantages of this paper is that we comprise the basic spline with classical polynomials both on accuracy, as well as degree of paralleling calculations are also shown.**

## I. INTRODUCTION

The mathematical device of splines developed within last decades by efforts of many researchers, has taken a worthy place among methods and algorithms of digital signal processing.

Splines as a class piece functions owing to universality of algorithms of processing of readout, good differential and extreme properties, high convergence of estimations, simplicity of calculations of forms and parameters, weak influence of mistakes of a rounding off find more and more wide application at creation of equipment rooms and software of the analysis and restoration of one-dimensional and multivariate signals, expanding frameworks of traditional approaches.

In comparison with approximation the classical polynomials the spline - function have two important advantages. At first, best approximate properties and, secondly, by convenience of realization of the algorithms, constructed on their basis, on the computer [1], [3], [4], [6].

## II. CUBIC BASIC SPLINES

Traditional methods of approximation of complex functional dependences have restrictions on accuracy, smooth nesses and do not allow essentially parallel computing structures for realization of approximation by classical polynomials on algorithm Horner [1], [4], [7], [8].

Cubic basic spline at $h=1$ are set by the following expressions, and the graph is represented in a fig.1:

$$B_3(x) = \begin{cases} 0, & x \geq 2 \\ (2-x)^3/6, & 1 \leq x \\ (1+3(1-x)+3(1-x)^2-3(1-x)^3)/6, & 0 \leq x \\ B_3(-x), & x < 0 \end{cases} \quad (1)$$

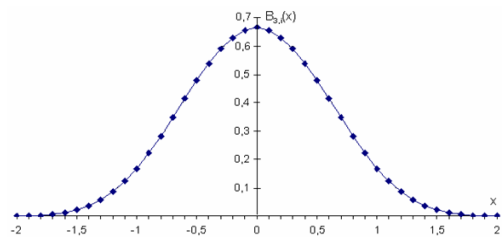

Fig. 1 The graph of a cubic basic spline

For an approximation of one dimensional function $f(x)$ with system of polynomial basic splines on an interval $[a,b]$ is necessary to enter outside an interval $[a,b]$ additional points in an amount 2m. Then the formula of an approximation of function $f(x)$ with basic splines can be written as the sum [1], [2], [5], and [6]:

$$f(x) \cong S_m(x) = \sum_{i=-1}^{m+1} b_i \cdot B_i(x), \quad a \leq x \leq b \quad (2)$$

$$f(x) \cong S_3(x) = b_{-1}B_{-1}(x) + b_0B_0(x) + b_1B_1(x) + b_2B_2(x) \quad (3)$$

Where: m – degree of a spline, $b_i$ - coefficients of smoothing, $B_i(x)$ - basic spline, $f(x)$ - function, $S_3(x)$ - cubic spline.

The three-dot formula, in this case for calculation of coefficient participate r-1, r and r+1-th values of function

$$b_r = \frac{1}{6}(-f_{r-1} + 8f_r - f_{r+1}) \qquad (4)$$

Here $f_r$ the reduced form of $f_r(x)$.

Thus, the local properties of basis is appears, and the structure can be created according to tabular - algorithmic methods.

The formula (3) follows, that it is necessary to summarize pair products by groups till four and periodically to update these sums in a summator.

The block scheme (fig.2) consists of a pre-set block, clock pulse generator (GEN), address counter (CT), generator of b-coefficients, shift register, memory (MEM) consisting of four subsections (ROM1, ROM2, ROM3, ROM4) for storage of value of a basic spline, four multipliers (M) and summator (SUM).

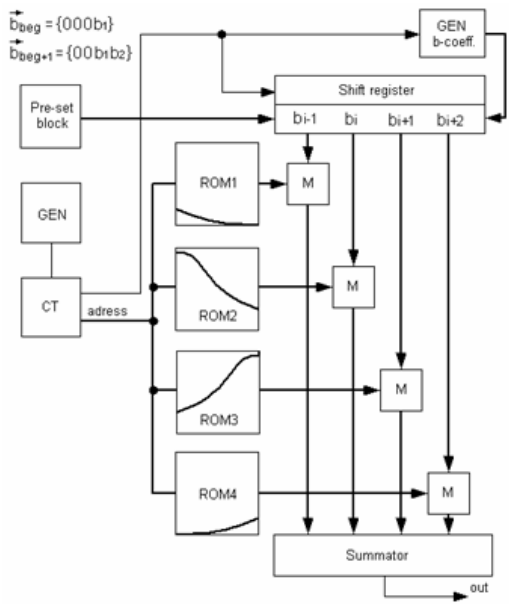

Fig. 2. The block scheme parallel specialized calculating machine on the basis of a cubic basic spline for restoring signals

All basic splines are identical. Hence in the storing device is possible to keep the table of values only of one basic spline. In fig 3 the cubic basic spline is represented, which is defined as nonzero on four segments of a piecewise partition. For samples of four different values of basic splines for summarize with the formula (3) will need four subsections of MEM, in each of subsection should be stored ¼ part of a curve of a basic spline given on one segment.

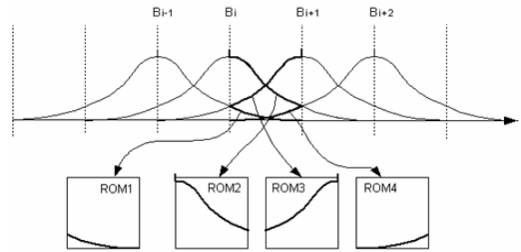

Fig. 3 The scheme of a memory allocation for appropriate segments of a basic spline.

If to keep in the MEM of values only of one basic spline, then with the two digits of a binary code is possible to define values of other basic splines. For sample of values of basic splines in each of four subsections of MEM, is possible to appropriate the following addresses:

ROM1 - 00 for an interval $x \in [1, 2)$;

ROM2 - 01 for an interval $x \in [0, 1)$;

ROM3 - 10 for an interval $x \in [-1, 0)$;

ROM4 - 11 for an interval $x \in [-2, -1)$.

The values of b-coefficient are worked out on the generator of b-coefficients with the formula (4) and records to the shift register accordingly by blocks $b_{i-1}, b_i, b_{i+1}, b_{i+2}$.

In the beginning the block of pre-set in the shift register establishes initial values of b-coefficients:

$$\vec{b}_{beg} = \{0\ 0\ 0\ b_1\}$$

$$\vec{b}_{beg+1} = \{0\ 0\ b_1\ b_2\}$$ etc.

Further in each cycle of operations the values of b-coefficients are established with shift one digit to left in the shift register.

In an initial step of the clock pulse generator the address counter specifies addresses of first values of a basic spline $B_i(x)$ kept in ROM1, ROM2, ROM3 and ROM4. Further in the following steps the CT specifies addresses of the following values of basic splines. After every ten cycles these values are iterate, since in our case four subsections of memory are used and in each memory ten values of a cubic basic spline are placed (provided that the exactitude of approximating is given by a constant).

Further appropriate values of basic splines and b-coefficients are incoming to multipliers (M) and parallel multiplies.

Hence, outcomes of multiplication are summarized in the summator and the data's out of summator are $S_3(x)$.

## III. Modeling of High-Efficiency Computing Structure

The modeling of the above mentioned High efficiency computing structure has been implemented using Simulink environment on MATLAB 7.4.0. A cubic spline was generated using the formula (1) and the b-coefficients were generated using formula (4). The four subsections of the spline were stored in the MATLAB workspace. The model is shown in fig. 4. A shift-left register was used to shift the b-coefficients at each iteration.

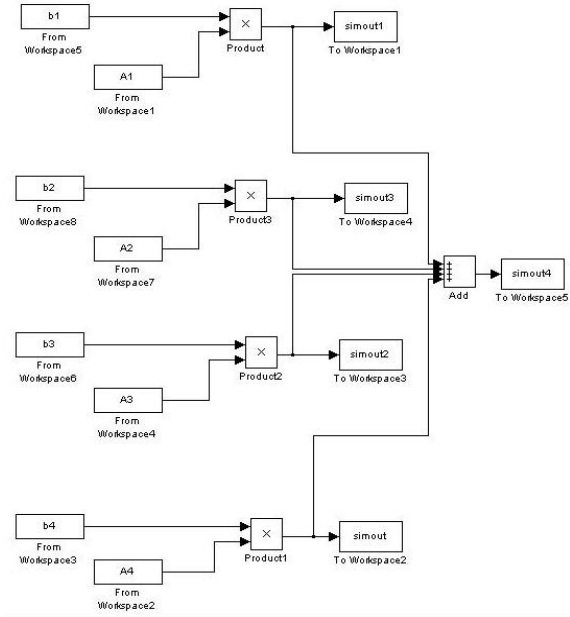

Fig. 4. The Simulink Model of the High efficiency computing structure

An implementation of the function $f(x) = \ln(1 + x)$ has been done using the structure and shown in fig. 5. The following section compares the function errors with the original function in the range [0, 2].

## IV. Comparison of the Received Results

The methodical error of interpolation of function $f(x)$ by cubic basic spline is defined by an inequality:

$$\varepsilon \le \frac{5}{384} h^4 \max \left| f^{IV}(x) \right| \qquad (5)$$

For function $f(x) = \ln(1 + x)$, let's receive

$$\varepsilon \le \frac{5}{384 \cdot 1{,}0 \cdot 32^4} = 0{,}12 \cdot 10^{-7} \qquad (6)$$

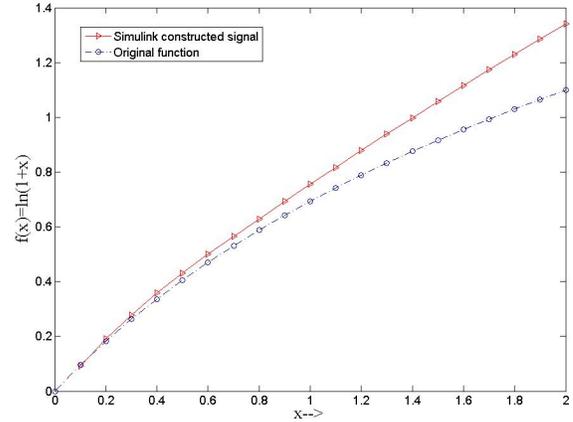

Fig. 5. Simulation result using the high efficiency computing structure for the function $f(x) = \ln(1 + x)$

For comparison we shall result value of an error of interpolation by classical cubic polynomials:

$$\varepsilon \le \frac{1}{24} h^4 \max \left| f^{IV}(x) \right| = \frac{1}{24 \cdot 32^4 \cdot 1{,}0} = 0{,}4 \cdot 10^{-7} \qquad (7)$$

Apparently from (6), the error exceeds the size received in (7), more than three times.

## V. Conclusion

Thus, (fig.3.) the high-efficiency computing structure is developed for restoration of the signals. A modeling of the structure using MATLAB simulink environment shows a high efficiency both in terms of the speed and the error. It functions three times more quickly than the existing, and the methodical error of the existing device exceeds three times, than an error of offered structure.

### Acknowledgment

This work was supported by the MIC (Ministry of Information and Communication), Korea, Under the ITFSIP (IT Foreign Specialist Inviting Program) supervised by the IITA (Institute of Information Technology Advancement)